\input amstex
\input epsf
\documentstyle{amsppt}

\def\B{\mathop{\text{\rm B}}}
\def\Im{\mathop{\text{\rm Im}}}
\def\Isom{\mathop{\text{\rm Isom}}}
\def\PU{\mathop{\text{\rm PU}}}
\def\Re{\mathop{\text{\rm Re}}}
\def\S{\mathop{\text{\rm S}}}
\def\SU{\mathop{\text{\rm SU}}}
\def\ta{\mathop{\text{\rm ta}}}
\def\tr{\mathop{\text{\rm tr}}}
\def\U{\mathop{\text{\rm U}}}

\hsize450pt

\topmatter
\title On the type of triangle groups
\endtitle\author Carlos H.~Grossi\endauthor
\thanks Supported by FAPEMIG.\endthanks
\address Departamento de Matem\'atica, ICEX, Universidade Federal de
Minas Gerais,\newline31161-970--Belo Horizonte--MG, Brasil\endaddress
\email grossi$_-$ferreira\@yahoo.com\endemail
\abstract We prove a conjecture of R.~Schwartz about the type of some
complex hyperbolic triangle groups.\endabstract
\endtopmatter\document

\centerline{\bf1.~Introduction}

\medskip

An $(n_1,n_2,n_3)$-complex reflection complex hyperbolic triangle group
is a group of isometries of the complex hyperbolic plane
$\Bbb H^2_\Bbb C$ generated by complex reflections $I_1,I_2,I_3$ in
complex geodesics $C_1,C_2,C_3$ such that $C_i$ and $C_{i+1}$ meet at
the angle $\pi/n_i$, $n_i\ge2$ (see Section 2 for definitions). For
fixed $n_1,n_2,n_3$, modulo conjugacy in $\Isom(\Bbb H^2_\Bbb C)$,
there exists in general a $1$-parameter family of
$(n_1,n_2,n_3)$-triangle groups. Assume $n_1\le n_2\le n_3$,
$n_i\in\Bbb N$. The triple $(n_1,n_2,n_3)$ is classified with respect
to the behavior of the isometries
$$W_A:=I_3I_2I_1I_2\quad\text{\rm{and}}\quad W_B:=I_1I_2I_3$$
while a parameter of the $(n_1,n_2,n_3)$-triangle group family varies
in a canonical way. The triple is said to be of type $A$ if $W_A$
becomes regular elliptic before $W_B$ and of type $B$ if $W_B$ becomes
regular elliptic before $W_A$ (see Subsections 2.2 and 2.3).

\smallskip

In this paper, we prove the following

\medskip

{\bf1.1.~Conjecture {\rm[Sch2, Conjecture 5.2]}.}~{\sl The triple\/
$(n_1,n_2,n_3)$ has type\/ $A$ if\/ $n_1\le9$ and type\/ $B$ if\/~
$n_1\ge14$.}

\medskip

Conjecture 1.1 is a tiny portion of a complete conjectural picture
[Sch2] that describes $(n_1,n_2,n_3)$-triangle groups with focus on
discreteness. Roughly speaking, the picture is as follows: An
$(n_1,n_2,n_3)$-triangle group is discrete if neither $W_A$ nor $W_B$
is regular elliptic. When $(n_1,n_2,n_3)$ has type $B$, then the
converse also holds. If $(n_1,n_2,n_3)$ has type $A$, then there is a
countable collection of `extra' discrete groups. While sorting triangle
groups into types $A$ and $B$ is (as we shall see) a matter of a couple
of simple tricks intended to avoid the huge amount of calculus a
straightforward approach leads to, classifying those groups with respect
to discreteness is certainly a much more difficult and interesting task.
It has been accomplished for $(\infty,\infty,\infty)$-triangle groups in
[GoP] and [Sch1] (see also [Sch3]) and for $(n_1,n_2,n_3)$-triangle
groups with sufficiently large $n_1$ in [Sch4].

\smallskip

In principle, we could say that our proof is computer independent. In
fact, we use the computer only to obtain approximate values of the
cosine function needed to establish some inequalities that hold by a
wide margin of error (Lemmas 3.1 and 3.2). The inequalities proved in
this way are marked with the symbols $\prec$ and $\succ$.

\smallskip

Conjecture 1.1 has been solved for sufficiently large $n_1$ in [Sch4]
and for triples of the form $(n,n,\infty)$ in [W-G]. In [Pra], it is
shown that the triples corresponding to triangles with
$r_1^2+r_2^2+r_3^2-1=2r_1r_2r_3$ and
$r_1r_2r_3\ge\frac{13+\sqrt{297}}{32}$ (see [Pra] and also Subsection
2.2) are of type $B$.

\medskip

{\bf Acknowledgments.} The author is grateful to Sasha Anan$'$in and
to Nikolay Gusevskii for many stimulating and fruitful discussions
about complex hyperbolic geometry.

\bigskip

\centerline{\bf2.~Preliminaries}

\medskip

{\bf2.1.~Basic background.}~Let $V$ be a $3$-dimensional $\Bbb
C$-vector space equipped with a hermitian form~
$\langle-,-\rangle$ of signature $+,+,-$. The complex hyperbolic
plane $\Bbb H^2_\Bbb C$ can be identified with the open $4$-ball
\footnote{The symbol := stands for `equals by definition.'}
$$\B V:=\big\{p\in\Bbb P_\Bbb CV\mid\langle p,p\rangle<0\big\}.$$
The complex hyperbolic distance $d(p_1,p_2)$ between two points
$p_1,p_2\in\B V$ is given in terms of the {\it tance\/}
\footnote{Here, and in what follows, we frequently do not
distinguish the notation of a point in $\Bbb P_\Bbb CV$ and of a
chosen representative of it in $V$ when a concept or expression does
not depend on such a choice.}
$$\ta(p_1,p_2):=\frac{\langle p_1,p_2\rangle\langle
p_2,p_1\rangle}{\langle p_1,p_1\rangle\langle p_2,p_2\rangle}$$ by
$\cosh^2\big(d(p_1,p_2)/2\big)=\ta(p_1,p_2)$ [Gol, p.~77]. The
ideal boundary of $\Bbb H^2_\Bbb C$ is the $3$-sphere
$$\S V:=\big\{p\in\Bbb P_\Bbb CV\mid\langle p,p\rangle=0\big\}$$
formed by the {\it isotropic\/} points in $\Bbb P_\Bbb CV$. Notice
that the tance $\ta(p_1,p_2)$ is well-defined for all
nonisotropic $p_1,p_2\in\Bbb P_\Bbb CV$. We put $\overline\B V:=\B
V\cup\S V$.

\smallskip

Every projective line $L$ in $\Bbb P_\Bbb CV$ has the form
$L=\Bbb P_\Bbb Cp^\perp$, where $p\in\Bbb P_\Bbb CV$ and
$p^\perp=\{v\in V\mid\allowmathbreak\langle v,p\rangle=0\}$. We call
$p$ the {\it polar point\/} to $L$. If $p\notin\overline\B V$, then
$\Bbb P_\Bbb Cp^\perp\cap\overline\B V$ is a {\it complex geodesic.}
Two distinct complex geodesics $C_1,C_2$ are {\it concurrent\/}
(respectively, {\it asymptotic,} {\it ultraparallel\/}) if and
only if $C_1\cap C_2\in\B V$ (respectively, $C_1\cap C_2\in\S V$,
$C_1\cap C_2=\varnothing$).

\medskip

{\bf2.1.1.~Lemma {\rm[Gol, p.~100]}.}~{\sl Two distinct complex
geodesics\/ $C_1,C_2$ with polar points\/ $p_1,p_2$ are concurrent,
asymptotic, ultraparallel if and only if\/ $\ta(p_1,p_2)<1$,
$\ta(p_1,p_2)=1$, $\ta(p_1,p_2)>1$, respectively. If\/
$\ta(p_1,p_2)\le1$, then the angle\/ $\angle(C_1,C_2)\in[0,\pi/2]$
between\/ $C_1$ and\/ $C_2$ is given by\/
$\cos^2\angle(C_1,C_2)=\ta(p_1,p_2)$}
$_\blacksquare$

\medskip

Given $p\notin\overline\B V$, define $I\in\SU V$ by the rule
$$I:x\mapsto2\frac{\langle x,p\rangle}{\langle p,p\rangle}p-x.
\leqno{\bold{(2.1.2)}}$$ The corresponding isometry in $\PU V$ is
known as the {\it complex reflection\/} in the complex geodesic
$\Bbb P_\Bbb Cp^\perp\cap\overline\B V$. For brevity, we will call
it simply the {\it reflection\/} in the complex geodesic in
question.

\medskip

{\bf2.2.~Complex hyperbolic triangles.}~A {\it complex hyperbolic
triangle\/} is a triple $(C_1,C_2,C_3)$ of complex geodesics in
$\overline\B V$. Each complex geodesic $C_i$ is a {\it side\/} of
the triangle. If the sides $C_{i}$ and $C_{i+1}$ meet at the angle
${\pi}/n_i$, where $n_i\ge2$ (we allow $n_i$ to be
infinite, meaning that $C_i$ and $C_{i+1}$ are asymptotic or
equal), the triangle $(C_1,C_2,C_3)$ is referred to as an
{\it{\rm(}$n_1,n_2,n_3${\rm)}-triangle.} We call an
$(n_1,n_2,n_3)$-triangle $(C_1,C_2,C_3)$ {\it
nondegenerate\/} if the form restricted to the subspace
$\Bbb Cp_1+\Bbb Cp_2+\Bbb Cp_3$ of $V$ is nondegenerate, being $p_i$
the polar point to $C_i$.

\medskip

{\bf2.2.1.~Lemma {\rm(compare with [Pra, Proposition 1])}.}~
{\sl Let\/ $(C_1,C_2,C_3)$ be a nondegenerate\/
$(n_1,n_2,n_3)$-triangle with\/ $n_i>2$,
$i=1,2,3$. Denote by\/ $p_i$ the polar point to\/ $C_i$. Define
$$r_i:=\sqrt{\ta(p_i,p_{i+1})}=\cos\frac{\pi}{n_i}>0,\quad
\varkappa:=\frac{\langle p_1,p_2\rangle\langle
p_2,p_3\rangle\langle p_3,p_1\rangle}{\langle
p_1,p_1\rangle\langle p_2,p_2\rangle\langle
p_3,p_3\rangle},\quad\varepsilon:=\frac{\varkappa}{|\varkappa|},\quad
t:=\Re\varepsilon.$$ Then, the numbers\/ $r_i$ and\/ $t$
constitute a complete set of geometrical invariants of\/
$(C_1,C_2,C_3)$. They satisfy\/ $0<r_i\le1$, $|t|\le1$, and
$$1+2r_1r_2r_3t-(r_1^2+r_2^2+r_3^2)\le0.\leqno{\bold{(2.2.2 )}}$$
All values of\/ $r_i$ and\/ $t$ subject to the conditions\/
$0<r_i\le1$, $|t|\le1$, and\/
$1+2r_1r_2r_3t-(r_1^2+r_2^2+r_3^2)\le0$ correspond to a
nondegenerate\/ $(n_1,n_2,n_3)$-triangle with\/
$n_i>2$.}

\medskip

{\bf Proof.}~The numbers $r_i$ and $\varepsilon$ are invariant
under the action of $\PU V$ on the triple $(p_1,p_2,p_3)$.
Choosing suitable representatives $p_i\in V$, we can assume that
$$\left(\matrix1&r_1&r_3\overline\varepsilon\\r_1&1&r_2\\r_3
\varepsilon&r_2&1\endmatrix\right)\leqno{\bold{(2.2.3)}}$$ is the
Gram matrix of $(p_1,p_2,p_3)$. If the triples $(p_1,p_2,p_3)$ and
$(p'_1,p'_2,p'_3)$ have the same Gram matrix and if the hermitian
form is nondegenerate being restricted to the subspaces generated
by $p_1,p_2,p_3$ and by $p'_1,p'_2,p'_3$, then there exists
$I\in\U V$ such that $Ip_i=p'_i$. The triangles corresponding to
$(r_1,r_2,r_3,\varepsilon)$ and to $(r'_1,r'_2,r'_3,\varepsilon')$
differ by an anti-holomorphic isometry of $\Bbb H^2_\Bbb C$ if and
only if $r_i=r'_i$ and $\varepsilon=\overline\varepsilon'$. The
rest follows from Sylvester's criterion $_\blacksquare$

\medskip

From now on, all $(n_1,n_2,n_3)$-triangles are assumed
to be nondegenerate.

\medskip

{\bf2.2.4.~Remark.}~For fixed $3<n_1\le n_2\le n_3$, there exists a
non-empty $1$-parameter family of $(n_1,n_2,n_3)$-triangles. Indeed,
the left-hand side of the inequality (2.2.2) is increasing in $t$ and,
hence, attains its minimum value at $t=-1$. We have
$$\displaystyle1-2r_1r_2r_3-(r_1^2+r_2^2+r_3^2)<1-2\cos^3\frac{\pi}
{3}-3\cos^2\frac{\pi}{3}=0\,_\blacksquare$$

Let $3<n_1\le n_2\le n_3$. In the terms of Lemma
2.2.1, define
$$t_M:=\frac{r_1^2+r_2^2+r_3^2-1}{2r_1r_2r_3},\quad
t_{max}:=\min\{t_M,1\}.\leqno{\bold{(2.2.5)}}$$
The {\it canonical\/} path of deformation of the
$(n_1,n_2,n_3)$-triangle family is the one that starts with $t=-1$ and
ends with $t=t_{max}$.

\medskip

{\bf2.3.~Complex hyperbolic triangle groups.}~The subgroup in $\PU V$
generated by the reflections in the sides of an
$(n_1,n_2,n_3)$-triangle is called an {\it$(n_1,n_2,n_3)$-triangle
group.} Up to conjugacy in $\Isom(\Bbb H^2_\Bbb C)$, all nondegenerate
$(n_1,n_2,n_3)$-triangle groups are described by Lemma 2.2.1.

\smallskip

Fix $3<n_1\le n_2\le n_3$ and assume $n_i\in\Bbb N$. For a given
$(n_1,n_2,n_3)$-triangle $(C_1,C_2,C_3)$, define
$$W_A:=I_3I_2I_1I_2,\quad W_B:=I_1I_2I_3,$$
where $I_i$ stands for the reflection in $C_i$. The triple
$(n_1,n_2,n_3)$ is characterized with respect to the behavior of the
isometries $W_A$ and $W_B$ during the canonical deformation of the
$1$-parameter family of $(n_1,n_2,n_3)$-triangle groups. Specifically,
$(n_1,n_2,n_3)$ is said to be of {\it type\/ $A$} if $W_A$ becomes
regular elliptic before $W_B$ and of {\it type\/ $B$} if $W_B$ becomes
regular elliptic before $W_A$ (see [Gol] for the classification of
holomorphic isometries of $\Bbb H^2_\Bbb C$).

\medskip

{\bf2.3.1.~Lemma {\rm(compare with [Pra, Proposition 12])}.}~{\sl
$W_A$ is always hyperbolic at the beginning of the deformation. In the
terms of Lemma\/~{\rm2.2.1}, $W_A$ is regular elliptic if and only if
$$t>t_{W_A}:=\frac{r_3^2+4r_1^2r_2^2-1}{4r_1r_2r_3}.$$}
\indent{\bf Proof.}~Let $(C_1,C_2,C_3)$ be an
$(n_1,n_2,n_3)$-triangle and let $p_i$ denote the
polar point to $C_i$. The isometry $W_A$ is the product of two
reflections: one in the complex geodesic with polar point $I_2p_1$
and the other in $C_3$. The nature of $W_A$ is hence determined by
the relative position of these complex geodesics.  Taking (2.2.3)
as the Gram matrix of suitable representatives $p_i\in V$ and
applying (2.1.2), we obtain
$$\ta(I_2p_1,p_3)=\frac{\langle I_2p_1,p_3\rangle\langle
p_3,I_2p_1\rangle}{\langle I_2p_1,I_2p_1\rangle\langle
p_3,p_3\rangle}=\big|\langle 2r_1p_2-p_1,p_3\rangle\big|^2=
4r_1^2r_2^2-4r_1r_2r_3t+r_3^2.$$ At the beginning of the
deformation,
$$\ta(I_2p_1,p_3)=4r_1^2r_2^2+4r_1r_2r_3+r_3^2>4(\cos^4\frac{\pi}{3}+
\cos^3\frac{\pi}{3})+\cos^2\frac{\pi}{3}=1.$$ By Lemma 2.1.1, this
implies that $W_A$ is hyperbolic. It remains to observe that
$\ta(I_2p_1,p_3)$ is decreasing in $t$ and that, by Lemma 2.1.1,
$W_A$ becomes parabolic exactly when $t=t_{W_A}$ $_\blacksquare$

\medskip

In order to deal with $W_B$, we need the following

\medskip

{\bf2.3.2.~Lemma {\rm[Gol, Theorem 6.2.4]}.}~{\sl Define a map\/
$f:\Bbb C\to\Bbb R$ by
$$f(z):=|z|^4-8\Re(z^3)+18|z|^2-27.$$
Given\/ $J\in\PU V$, let\/ $\hat J\in\SU V$ be a lift of\/ $J$.
Then, $J$ is regular elliptic\/ {\rm(}respectively,
loxodromic\/{\rm)} if and only if\/ $f(\tr\hat J)<0$
{\rm(}respectively, $f(\tr\hat J)>0${\rm).} The pre-image\/
$f^{-1}(0)\subset\Bbb C$ is a deltoid. If\/ $\tr\hat J\in\Bbb R$,
then\/ $J$ is loxodromic if and only if\/
$\tr\hat J\notin[-1,3]$
$_\blacksquare$}

\medskip

Take the lift $W_B\in\SU V$ determined by the lifts of
$I_1,I_2,I_3$ in (2.1.2). The trace $\tau:=\tr W_B$ is given by
$$\tau:=\tr W_B=8r_1r_2r_3\varepsilon-4(r_1^2+r_2^2+r_3^2)+3\in\Bbb C
\leqno{\bold{(2.3.3)}}$$
(see, for instance, [Pra]).

\medskip

{\bf2.3.4. Lemma.}~{\sl$W_B$ is always loxodromic at the beginning of
the deformation.}

\medskip

{\bf Proof.} By (2.3.3),
$$\tau=-8r_1r_2r_3-4(r_1^2+r_2^2+r_3^2)+3<-8\cos^3\frac{\pi}{3}-12
\cos^2\frac{\pi}{3}+3=-1$$
at the beginning of the deformation. The result follows from
Lemma 2.3.2
$_\blacksquare$

\medskip

During the deformation, $\tau$ belongs to the circle
$$F:=\Big\{(x,y)\in\Bbb C\mid\big(x+4(r_1^2+r_2^2+r_3^2)-3\big)^2+
y^2=(8r_1r_2r_3)^2\Big\},\leqno{\bold{(2.3.5)}}$$ where
$x:=\Re\tau$ and $y:=\Im\tau$. By Lemma 2.2.1, we can assume that
$\Im\varepsilon\ge0$, i.e., $y\ge0$. By (2.3.3), the coordinate
$x=\Re\tau$ and the parameter $t=\Re\varepsilon$ are linked by the
formula
$$t=\frac{x+4(r_1^2+r_2^2+r_3^2)-3}{8r_1r_2r_3}.
\leqno{\bold{(2.3.6)}}$$
Hence, we can think of the parameter $t$ as `living' in the upper
half-circle of $F$.

\centerline{\bf3.~Proof of the conjecture}

\medskip

In what follows, we will refer to the elementary Lemmas 4.1 and 4.3,
proved in Section 4.

\medskip

{\bf3.1.~Proposition.}~{\sl If\/ $14\le n_1\le n_2\le n_3$, then the
triple\/ $(n_1,n_2,n_3)$ is of type\/ $B$.}

\medskip

{\bf Proof.}~In the terms of Lemma 2.2.1,
$\cos\frac{\pi}{14}\le r_1\le r_2\le r_3\le1$.

\smallskip

Intersection points $(x,-\sqrt3x)\in\Bbb C$ of the line $y=-\sqrt3x$
that passes through the vertex $(-\frac{3}{2},\frac{3\sqrt3}{2})$ and
through the center $(0,0)$ of Goldman's deltoid (Lemma 2.3.2) with the
circle $F$ given by (2.3.5) satisfy\footnote{Obviously, we always sum
over $i=1,2,3$.}
$$4x^2+2(4\sum r_i^2-3)x+\big(4\sum r_i^2-3\big)^2-(8r_1r_2r_3)^2=0.$$
By Lemma 4.1 (1), the discriminant
$$D_1:=4\big(-3(4\sum r_i^2-3)^2+(16r_1r_2r_3)^2\big)$$
of the above equation is such that $D_1>0$. Take the root
$\displaystyle x_0:=\frac{-2(4\sum r_i^2-3)+\sqrt{D_1}}{8}.$

\medskip

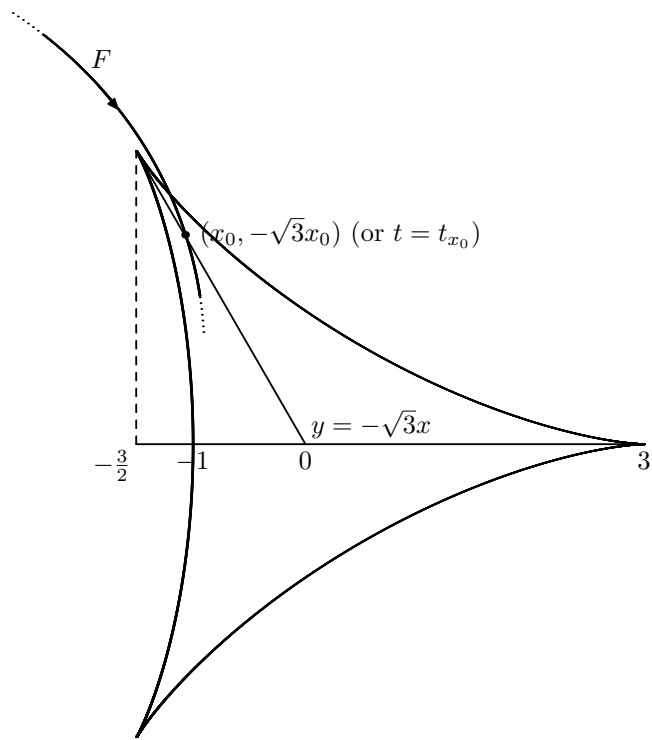
\centerline{\noindent$\vcenter{\hbox{\epsfbox{deltoid.1}}}$}

\medskip

We will show that $x_0\in(-\frac{3}{2},-1)$. In particular, this
implies that $(x_0,-\sqrt3x_0)\in F$ is `inside of' Goldman's deltoid.
\footnote{It would suffice to prove here a weaker inequality, but we
will need later the fact that $x_0<-1$.} The inequality $x_0<-1$,
equivalent to $\sqrt{D_1}<2(4\sum r_i^2-3)-8$, follows from
Lemma~4.1~(2) and from\footnote{As stated in the introduction, the
symbols $\prec$ and $\succ$ are used for the inequalities proved using
the computer to find approximate values of the cosine function.}
$$2(4\sum r_i^2-3)-8\ge2(12\cos^2\frac{\pi}{14}-3)-8\succ8.8.$$
The inequality $x_0>-\frac{3}{2}$ follows from
$-(4\sum r_i^2-3)^2+3(4\sum r_i^2-3)+(8r_1r_2r_3)^2-9>0$ which is a
consequence of Lemma 4.1 (3).

\smallskip

According to (2.3.6), the value of the deformation parameter $t$ that
corresponds to $x_0$ is
$$t_{x_0}:=\frac{x_0+4\sum_ir_i^2-3}{8r_1r_2r_3}.$$
It satisfies $t_{x_0}\in(-1,t_{max})$, being $t_{max}$ as defined in
(2.2.5). Indeed, the inequality $t_{x_0}>-1$ is straightforward. The
inequality $t_{x_0}<1$, equivalent to
$\sqrt{D_1}<-6\big(4\sum r_i^2-3\big)+64r_1r_2r_3$, follows from
Lemma~4.1~(2),~(4). Finally, the inequality $t_{x_0}<t_M$ is
equivalent to $x_0<-1$.

\smallskip

We have just proved that the deformation parameter assumes the value
$t=t_{x_0}$. By Lemma 2.3.2, the isometry $W_B$ is regular elliptic
when $t=t_{x_0}$ since $(x_0,-\sqrt3x_0)$ is inside of the deltoid. By
Lemma~2.3.1, in order to show that $W_B$ becomes regular elliptic
before $W_A$, it suffices to show that $t_{x_0}<t_{W_A}$. This follows
from
$$3(4\sum r_i^2-3)^2-12(r_3^2+4r_1^2r_2^2-1)(4\sum r_i^2-3)+
16(r_3^2+4r_1^2r_2^2-1)^2-(8r_1r_2r_3)^2>0$$
which is a consequence of Lemma 4.1 (5)
$_\blacksquare$

\medskip

{\bf3.5.~Proposition.}~{\sl If\/ $n_1\le n_2\le n_3$ and\/
$4\le n_1\le9$, then the triple\/ $(n_1,n_2,n_3)$ is of type\/ $A$.}

\medskip

{\bf Proof.} In the terms of Lemma 2.2.1, $r_1\le r_2\le r_3\le1$
and $\cos\frac{\pi}{4}\le r_1\le\cos\frac{\pi}{9}$.

\smallskip

First, let us show that the deformation parameter $t$ assumes the value
$t=t_{W_A}$, being $t_{W_A}$ as in Lemma~2.3.1. In other words, we need
to show that $t_{W_A}\in[-1,t_{max}]$ (see 2.2.5). The inequalities
$t_{W_A}>-1$ and $t_{W_A}\le1$ are straightforward. The inequality
$t_{W_A}\le t_M$, equivalent to~$2(r_1^2+r_2^2-2r_1^2r_2^2)+r_3^2-1\ge0$, is a consequence of Lemma
4.3 (1).

\medskip

Assume $\cos\frac{\pi}{4}\le r_1\le\cos\frac{\pi}{8}$.

\medskip

At the beginning $t=-1$ of the deformation, the trace $\tau$ of $W_B$
given by (2.3.3) satisfies
$$\tau\le-8\cos^3\frac{\pi}{4}-12\cos^2\frac{\pi}{4}+3<-\frac{3}{2}.$$
By Lemma 2.3.2, this means that $W_B$ may become elliptic only after
the parameter
$$t_{B}:=\frac{-\frac{3}{2}+4\sum r_i^2-3}{8r_1r_2r_3}$$
that corresponds, by (2.3.6), to $x=-\frac{3}{2}$.

\medskip

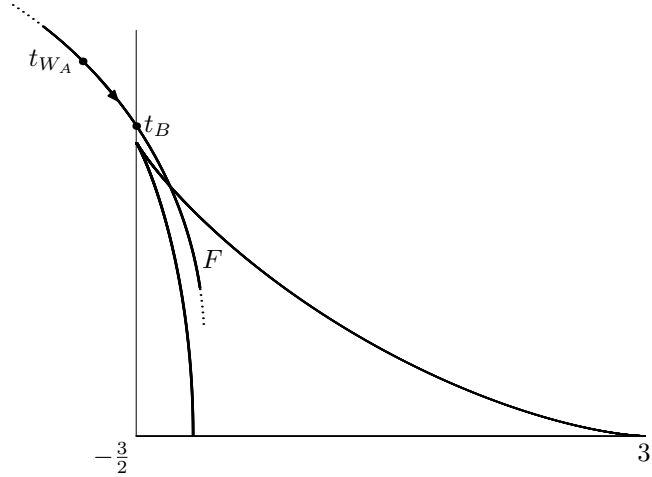
\centerline{\noindent$\vcenter{\hbox{\epsfbox{deltoid.2}}}$}

\medskip

\noindent The inequality $t_{W_A}<t_{B}$ is equivalent to
$8(r_1^2+r_2^2-2r_1^2r_2^2)+4r_3^2-5>0$ and follows from Lemma~4.3~
(2). This implies that $(n_1,n_2,n_3)$ is of type $A$.

\medskip

We now consider the case $r_1=\cos\frac{\pi}{9}$.

\medskip

By Lemma 2.3.2, intersection points of the deltoid with the line
$l:=\big\{(x,y)\in\Bbb C\mid y=\frac{3\sqrt3}{5}(1-x)\big\}$
satisfy $(2x+3)^2(169x^2-158x-111)=0$. The roots of this equation are
$x=-\frac{3}{2}$ (that corresponds
to a vertex of the deltoid),
$$x=x_1:=\frac{1}{169}(79-50\sqrt{10}),\quad\text{\rm and}
\quad x=\frac{1}{169}(79+50\sqrt{10}).$$

Intersection points of the circle $F$ given by (2.3.5) with the line
$l$ satisfy
$$52x^2+2\big(25(4\sum r_i^2-3)-27\big)x+25\big((4\sum r_i^2-3)^2-
(8r_1r_2r_3)^2\big)+27=0.$$
By Lemma 4.3 (3), the discriminant
$$D_2:=100\big(-27(4\sum r_i^2-3)^2-54(4\sum r_i^2-3)+
52(8r_1r_2r_3)^2-27\big)$$
of the above equation is such that $D_2>0$. Take the
root
$$x_2:=\frac{-2\big(25(4\sum r_i^2-3)-27\big)-\sqrt{D_2}}
{104}.$$

In order to prove that $(9,n_2,n_3)$ is of type $A$, it
suffices to apply Lemma 2.3.1 after showing the following facts (see
the picture below):

\medskip

(1) $x_2<-\frac{3}{2}$. This implies that, when $F$ crosses $l$ for
the first time (thus entering the region in grey), $W_B$ has not
become elliptic yet.

\smallskip

(2) $x_{W_A}<x_1$ and $g(x_{W_A},y_{W_A})>0$, where
$(x_{W_A},y_{W_A})\in F$ is the point that corresponds to $t_{W_A}$
by~(2.3.6) and $g(x,y):=y-\frac{3\sqrt3}{5}(1-x)$. This implies that
we are still in the grey region when $t=t_{W_A}$. In particular, $W_B$
has not become elliptic yet.

\bigskip

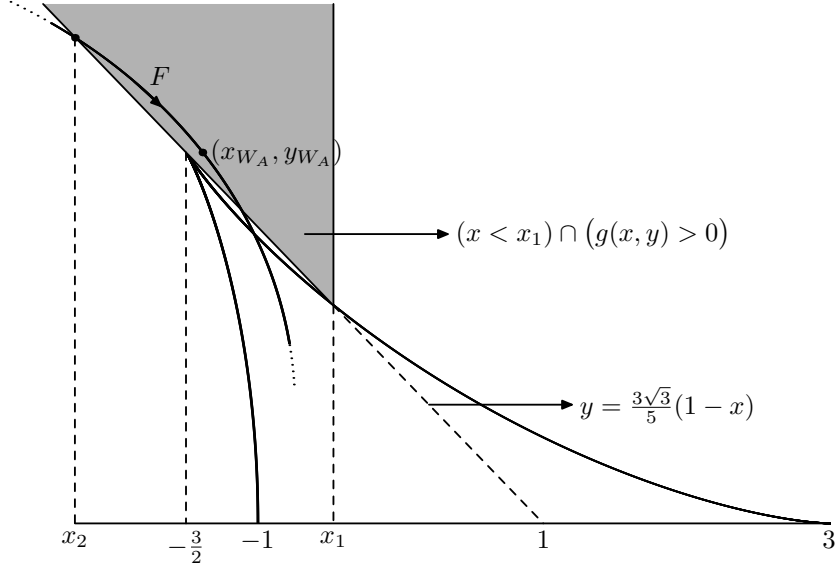
\centerline{\noindent$\vcenter{\hbox{\epsfbox{deltoid.3}}}$}

\bigskip

The inequality $x_2<-\frac{3}{2}$, equivalent to
$\sqrt{D_2}>-50(4\sum r_i^2-3)+210$, follows from
$$-50(4\sum r_i^2-3)+210\le-50(12\cos^2\frac{\pi}{9}-3)+210\prec0.$$

By (2.3.6) and (2.3.5),
$$x_{W_A}:=2(r_3^2+4r_1^2r_2^2-1)-(4\sum r_i^2-3),\quad
y_{W_A}:=2\sqrt{(4r_1r_2r_3)^2-(r_3^2+4r_1^2r_2^2-1)^2}.$$
The inequality $x_{W_A}<x_1$, equivalent to
$2(r_3^2+4r_1^2r_2^2-1)-(4\sum r_i^2-3)-\frac{1}{169}
(79-50\sqrt{10})<0$, follows from Lemma 4.3 (4). Finally,
$g(x_{W_A},y_{W_A})>0$ is a consequence of
$$25\big((4r_1r_2r_3)^2-(r_3^2+4r_1^2r_2^2-1)^2\big)-27
\big(2r_2^2(1-2r_1^2)+2r_1^2+r_3^2\big)^2>0$$
which follows from Lemma 4.3 (5)
$_\blacksquare$

\medskip

{\bf3.6.~Remark.}~For the sake of generality, we have not considered
yet the $(3,n_2,n_3)$-triangles, where $3\le n_2\le n_3\le\infty$. As
in Remark 2.2.4, it is easy to see that there exists in general a
non-empty one parameter family of $(3,n_2,n_3)$-triangles. The only
exception is the $(3,3,3)$-triangle, which is rigid. Proceeding as in
the first part of the proof of Proposition 3.5 (where we dealt with the
$4\le n_1\le 8$ cases) and using the fact that $n_i\in\Bbb N$, one
easily shows that the non-rigid $(3,n_1,n_2)$-triangles are of type $A$
$_\blacksquare$

\bigskip

\centerline{\bf4.~Taking Derivatives}

\bigskip

{\bf4.1.~Lemma.}~{\sl Suppose that\/
$\cos\frac{\pi}{14}\le x,y,z\le1$.
Define
\medskip
\noindent $f_1(x,y,z):=-3\big(4(x^2+y^2+z^2)-3\big)^2+(16xyz)^2,$
\smallskip
\noindent
$f_2(x,y,z):=-\big(4(x^2+y^2+z^2)-3\big)^2+3\big(4(x^2+y^2+z^2)-3
\big)+(8xyz)^2-9,$
\smallskip
\noindent $f_3(x,y,z):=-3\big(4(x^2+y^2+z^2)-3\big)+32xyz,$
\smallskip
\noindent
$f_4(x,y,z):=3\big(4(x^2+y^2+z^2)-3\big)^2-12(z^2+4x^2y^2-1)
\big(4(x^2+y^2+z^2)-3\big)+16(z^2+4x^2y^2-1)^2-(8xyz)^2$.
\medskip
Then,
$$\text{\bf(1)\;\;}31\prec4f_1(x,y,z)\le52,\quad
\text{\bf(2)\;\;}5.5<\sqrt{4f_1(x,y,z)}<7.3,\quad
\text{\bf(3)\;\;}f_2(x,y,z)\succ0.5,$$
$$\text{\bf(4)\;\;}2f_3(x,y,z)\succ8.8,\quad{\text{\rm and}}\quad
\text{\bf(5)\;\;}f_4(x,y,z)\succ0.1.$$}

{\bf Proof} is straightforward. We show, for instance, the first and
the last items. Notice that
$$8\cos^2\frac{\pi}{14}-3\succ0.\leqno{\bold{(4.2)}}$$
\vskip-2pt
\indent(1) Taking derivatives,
$$\frac{\partial f_1}{\partial x}=16x\Big(-3\big(4(x^2+y^2+z^2)-
3\big)+32y^2z^2\Big),\quad\frac{\partial^2f_1}{\partial x^2}=
16\Big(-36x^2-3\big(4(y^2+z^2)-3\big)+32y^2z^2\Big).$$
By (4.2),
$\displaystyle\frac{\partial^2f_1}{\partial x^2}$ is increasing in
$y$ and in $z$. Hence,
$\displaystyle\frac{\partial^2f_1}{\partial x^2}\le
\frac{\partial^2f_1}{\partial x^2}(\cos\frac{\pi}{14},1,1)\prec0$.
This implies that
$\displaystyle\frac{\partial f_1}{\partial x}$ is decreasing in $x$
and (4.2) implies that it is increasing in $y$ and in $z$. So,
$\displaystyle\frac{\partial f_1}{\partial x}\ge\displaystyle
\frac{\partial f_1}{\partial x}(1,\cos\frac{\pi}{14},\cos
\frac{\pi}{14})\succ0$.
We have just proved that $f_1$ is increasing in every variable. It
follows that
$$31\prec4f_1(\cos\frac{\pi}{14},\cos\frac{\pi}{14},\cos\frac{\pi}{14})
\le4f_1\le4f_1(1,1,1)\le52.$$

(5) Taking derivatives,
$$\frac{\partial f_4}{\partial x}=16x\big(4x^2(16y^4-12y^2+3)-2y^2
(12y^2+8z^2-7)+6z^2-3\big),$$
$$\frac{\partial^2f_4}{\partial x\partial
y}=64xy\big(-24(x^2+y^2)+64x^2y^2-8z^2+7\big),\quad\frac{\partial
f_4}{\partial z}=8z\big(12(x^2+y^2)-32x^2y^2+8z^2-5\big).$$
Put $g(y):=16y^4-12y^2+3$. Notice that $g'(y)>0\iff8y^2-3>0$ and
that the last inequality follows from (4.2). Hence,
$g(y)\ge g(\cos\frac{\pi}{14})\succ0$. This implies that
$\displaystyle\frac{\partial f_4}{\partial x}$ is increasing in
$x$. By (4.2),
$\displaystyle\frac{\partial^2f_4}{\partial x\partial y}$ is
increasing in $x$ and in $y$. Hence,
$\displaystyle\frac{\partial^2f_4}{\partial x\partial
y}\ge\frac{\partial^2f_4}{\partial x\partial
y}(\cos\frac{\pi}{14},\cos\frac{\pi}{14},1)\succ0.$
In other words, $\displaystyle\frac{\partial f_4}{\partial x}$ is
increasing also in $y$. It is decreasing in $z$ by (4.2). So,
$\displaystyle\frac{\partial f_4}{\partial x}\ge \frac{\partial
f_4}{\partial x}(\cos\frac{\pi}{14}, \cos\frac{\pi}{14},1)\succ0.$
It follows that $f_4$ is increasing in both $x$ and $y$. Moreover,
$\displaystyle\frac{\partial f_4}{\partial z}$ is increasing in
$z$ and decreasing in $x$ and in $y$ by~(4.2). This implies that
$\displaystyle\frac{\partial f_4}{\partial z}\le \frac{\partial
f_4}{\partial z}(\cos\frac{\pi}{14}, \cos\frac{\pi}{14},1)\prec0$,
that is, $f_4$ is decreasing in $z$. Finally,
$$f_4\ge f_4(\cos\frac{\pi}{14},\cos\frac{\pi}{14},1)\succ0.1
\,_\blacksquare$$

{\bf4.3.~Lemma.}~{\sl Define
\medskip
\noindent
$g_1(x,y,z):=2(x^2+y^2-2x^2y^2)+z^2-1$ for\/
$\cos\frac{\pi}{4}\le x\le y\le z\le1$,
\smallskip
\noindent
$g_2(x,y,z):=8(x^2+y^2-2x^2y^2)+4z^2-5$ for\/
$\cos\frac{\pi}{4}\le x\le\cos\frac{\pi}{8}$ and\/ $x\le y\le z\le1$.

\medskip

Suppose that\/ $x=\cos\frac{\pi}{9}$ and that\/
$\cos\frac{\pi}{9}\le y\le z\le1$. Define\/
\medskip
\noindent
$g_3(y,z):=-27\big(4(x^2+y^2+z^2)-3\big)^2-54\big(4(x^2+y^2+z^2)
-3\big)+52(8xyz)^2-27$,
\smallskip
\noindent
$g_4(y,z):=2(z^2+4x^2y^2-1)-\big(4(x^2+y^2+z^2)-3\big)-\frac{1}{169}
(79-50\sqrt{10})$,
\smallskip
\noindent
$g_5(y,z):=25\big((4xyz)^2-(z^2+4x^2y^2-1)^2\big)-27
\big(2y^2(1-2x^2)+2x^2+z^2\big)^2$.
\medskip
Then,
$$\text{\bf(1)\;\;}g_1(x,y,z)\ge0,\quad
\text{\bf(2)\;\;}g_2(x,y,z)\succ0.1,\quad
\text{\bf(3)\;\;}g_3(y,z)\succ296$$
$$\text{\bf(4)\;\;}g_4(y,z)\prec-0.9,\quad\text{and}\quad
\text{\bf(5)\;\;}g_5(y,z)\succ0.2.$$}
\noindent{\bf Proof} is straightforward. We show the first and the
last items, for instance.

\medskip

(1) Clearly, $g_1(x,y,z)\ge g_1(x,y,y)$. Since
$g_1(x,y,y)=2x^2(1-2y^2)+3y^2-1$ and
$1-2y^2\le0$, we obtain
$g_1(x,y,z)\ge g_1(x,y,y)\ge g_1(y,y,y)=-4y^4+5y^2-1\ge0$.

\smallskip

(5) Taking derivatives,
$$\frac{\partial g_5}{\partial y}=8y\Big(z^2(104x^2-27)+
2y^2\big(-100x^4-27(1-2x^2)^2\big)-4x^2(1-27x^2)\Big).$$
It follows from $104x^2-27\succ0$ and $-100x^4-27(1-2x^2)^2<0$ that
$\displaystyle\frac{\partial g_5}{\partial y}\le\frac{\partial
g_5}{\partial y}(\cos\frac{\pi}{9},1)\prec0$. In other words,
$g_5(y,z)$ is decreasing in $y$. Hence, $g_5(y,z)\ge g_5(z,z)$. Define
$$h(z):=g_5(z,z)=25\big(16x^2z^4-(z^2+4x^2z^2-1)^2\big)-
27\big(z^2(3-4x^2)+2x^2\big)^2.$$
We have
$$\frac{\partial h}{\partial z}=
4z\big(4z^2(-208x^4+212x^2-67)+25(1+4x^2)-54x^2(3-4x^2)\big).$$
It follows from
$-208x^4+212x^2-67\prec0$ that
$\displaystyle\frac{\partial h}{\partial
z}\le\displaystyle\frac{\partial h}{\partial z}(\cos\frac{\pi}{9})
\prec0$. So, $g_5(y,z)\ge h(z)\ge h(1)\succ0.2$ $_\blacksquare$

\bigskip

\centerline{\bf5.~References}

\medskip

[Gol] W.~M.~Goldman, {\it Complex Hyperbolic Geometry,} Oxford
Mathematical Monographs. Oxford Science Publications. The
Clarendon Press, Oxford University Press, New York, 1999

[GoP] W.~M.~Goldman, J.~R.~Parker, {\it Complex hyperbolic ideal
triangle groups,} J.~Reine Angew.~Math {\bf 425} (1992), 71--86

[Pra] A.~Pratoussevitch, {\it Traces in complex hyperbolic triangle
groups,} Geom.~Dedicata {\bf 111} (2005) 159--185

[Sch1] R.~E.~Schwartz, {\it Ideal triangle groups, dented tori, and
numerical analysis,} Ann.~of Math.~(2) {\bf 153} (2001), 533--598

[Sch2] R.~E.~Schwartz, {\it Complex Hyperbolic Triangle Groups,}
Proceedings of the International Congress of Mathematicians, Vol.~II,
Higher Ed.~Press, Beijing (2002) 339--349

[Sch3] R.~E.~Schwartz, {\it A better proof of the Goldman-Parker
conjecture,} Geom.~Topol.~{\bf9} (2005) 1539--1601

[Sch4] R.~E.~Schwartz, {\it Spherical\/ $CR$ Geometry and Dehn
Surgery,} Ann.~of Math.~Stud. {\bf165} 2007

[W-G] J.~Wyss-Gallifent, {\it Complex Hyperbolic Triangle Groups,}
Ph.~D.~thesis, University of Maryland (2000)

\enddocument